\documentclass[12pt ]{article}
\usepackage{amsthm,amsfonts,amsmath,amssymb}

\newcommand{\const}{\mathrm{const}}
\newcommand{\rd}{{\rm d}}
\newcommand{\re}{{\rm e}}

\newcommand{\myp}{\mbox{$\:\!$}}
\newcommand{\mypp}{\mbox{$\;\!$}}
\newcommand{\myn}{\mbox{$\:\!\!$}}

\newcommand{\RR}{{\mathbb R}}
\newcommand{\CC}{{\mathbb C}}

\newcommand{\PP}{{\mathsf P}}
\newcommand{\EE}{{\mathsf E}}

\newcommand{\calD}{{\mathcal D}}
\newcommand{\calF}{{\mathcal F}}

\newcommand{\calL}{{\mathcal L}}

\newtheorem{theorem}{Theorem}%[section]
\newtheorem{proposition}%[theorem]
{Proposition}
\newtheorem{definition}%[theorem]
{Definition}
\newtheorem*{corollary}{Corollary}

\makeatletter 

\newif\if@longtoc
\newif\if@indentheadings
\def\@seccntformat#1{\csname pre#1\endcsname\csname the#1\endcsname
                \csname post#1\endcsname}
\let\@Tocseccntformat\@seccntformat
\def\@postskip@{\hskip.5em\relax}
\def\postsection{.\@postskip@}
\def\postsubsection{.\@postskip@}
\def\postsubsubsection{.\@postskip@}
\def\postparagraph{.\@postskip@}
\def\postsubparagraph{.\@postskip@}

\def\@sect#1#2#3#4#5#6[#7]#8{%
     \ifnum #2>\c@secnumdepth
       \let\@svsec\@empty\else
       \refstepcounter{#1}%
       \let\@@protect\protect
       \def\protect{\noexpand\protect\noexpand}%
       \edef\@svsec{\@seccntformat{#1}}%
       \let\protect\@@protect
     \fi
     \@tempskipa #5\relax
      \ifdim \@tempskipa>\z@
        \begingroup #6\relax
          \@hangfrom{\hskip #3\relax\@svsec}%
                    {\interlinepenalty \@M #8\par}%
        \endgroup
       \csname #1mark\endcsname{#7}\addcontentsline
         {toc}{#1}{\ifnum #2>\c@secnumdepth \else
                      \protect\numberline{\@Tocseccntformat{#1}}\fi
         \if@longtoc#8\else#7\fi}\else
        \def\@svsechd{#6\hskip #3\relax  %% \relax added 2 May 90
                   \@svsec #8\csname #1mark\endcsname
                      {#7}\addcontentsline
                           {toc}{#1}{\ifnum #2>\c@secnumdepth \else
                           \protect\numberline{\@Tocseccntformat{#1}}\fi
         \if@longtoc#8\else#7\fi}}\fi
     \@xsect{#5}}

  \if@twoside
   \def\ps@headings{\let\@mkboth\markboth
   \let\@oddfoot\@empty\let\@evenfoot\@empty
   \def\@evenhead{\thepage\hfil\slshape\leftmark}%
   \def\@oddhead{{\slshape\rightmark}\hfil\thepage}%
    \def\sectionmark####1{\markboth {\uppercase{\ifnum \c@secnumdepth >\z@
      \thesection\postsection \hskip 1em\relax \fi ####1}}{}}%
    \def\subsectionmark####1{\markright {\ifnum \c@secnumdepth >\@ne
            \thesubsection\postsubsection \hskip 1em\relax \fi ####1}}}
  \else
   \def\ps@headings{\let\@mkboth\markboth
    \let\@oddfoot\@empty
    \def\@oddhead{{\slshape\rightmark}\hfil\thepage}%
    \def\sectionmark####1{\markright {\uppercase{\ifnum \c@secnumdepth >\z@
      \thesection\postsection \hskip 1em\relax \fi ####1}}}}
  \fi

\if@indentheadings
\def\section{\@startsection {section}{1}{\parindent}%
                                   {-3.5ex \@plus 1ex \@minus .2ex}%
                                   {2.3ex \@plus.2ex}%
                                   {\reset@font\large
                                   \bfseries}}
\def\subsection{\@startsection{subsection}{2}{\parindent}%
                                     {3.25ex\@plus 1ex \@minus .2ex}%
                                     {1.5ex \@plus .2ex}%
                                     {\reset@font%\large
                                     \itshape}}
\def\subsubsection{\@startsection{subsubsection}{3}{\parindent}%
                                     {3.25ex\@plus 1ex \@minus .2ex}%
                                     {1.5ex \@plus .2ex}%
                                     {\reset@font\normalsize\itshape}}
\else
\def\section{\@startsection {section}{1}{0pt}%
                                   {-3.5ex \@plus -1ex \@minus -.2ex}%
                                   {2.3ex \@plus .2ex}%
                                   {\reset@font\large
                                   \bfseries}}
\def\subsection{\@startsection{subsection}{2}{0pt}%
                                     {3.25ex\@plus -.5ex \@minus -.2ex}%
                                     {1.5ex \@plus .2ex}%
                                     {\reset@font%\large
                     \bfseries}}
\def\subsubsection{\@startsection{subsubsection}{3}{0pt}%
                                     {-1.5ex \@plus -.2ex\@minus -.2ex}%
                                     {-1.5ex \@plus -.2ex}%
                                     {\reset@font\normalsize\itshape}}
\fi

\makeatother

\begin{document}

\title{On Bounded Solutions of the Balanced Generalized Pantograph Equation}
\author{Leonid Bogachev\myp$^{\rm a}$, Gregory Derfel\myp$^{\rm b}$, Stanislav Molchanov\myp$^{\rm
c}$,
and\\ John Ockendon\myp$^{\rm d}$}
\date{\small $^{\,\rm a}$\myp Department of Statistics, University of
Leeds, Woodhouse Lane, Leeds LS2 9JT, UK.\\ E-mail: bogachev@maths.leeds.ac.uk\\
$^{\,\rm b}$\myp Department of Mathematics, Ben Gurion University of
the Negev, Beer Sheva, Israel. \\E-mail: derfel@math.bgu.ac.il\\
$^{\,\rm c}$\myp Department of Mathematics, University of North
Carolina at Charlotte, Charlotte\\ NC 28223, USA. E-mail: smolchan@uncc.edu\\
$^{\,\rm d}$\myp University of Oxford, Centre for Industrial and
Applied Mathematics, Mathematical Institute, 24-29 St Giles, Oxford
OX1 3LB, UK. E-mail: ock@maths.ox.ac.uk} \maketitle

\begin{abstract}
The question about the existence and characterization of bounded
solutions to linear functional-differential equations with both
advanced and delayed arguments was posed in early 1970s by T.~Kato
in connection with the analysis of the pantograph equation,
$y'(x)=a\myp y(q\myp x)+b\myp y(x)$. In the present paper, we answer
this question for the \emph{balanced} generalized pantograph
equation of the form $-a_2\myp y''(x)+a_1 y'(x)+y(x)=\int_0^\infty
y(\alpha x)\,\mu(\rd \alpha)$, where $a_1\ge0$, $a_2\ge0$,
$a_1^2+a_2^2>0$, and $\mu$ is a probability measure. Namely, setting
$K:=\int_0^\infty \ln \alpha\:\mu(\rd\alpha)$, we prove that if
$K\le 0$ then the equation does not have nontrivial (i.e.,
nonconstant) bounded solutions, while if $K>0$ then such a solution
exists. The result in the critical case, $K=0$, settles a
long-standing problem. The proof exploits the link with the theory
of Markov processes, in that any solution of the balanced pantograph
equation is an $\calL$-harmonic function relative to the generator
$\calL$ of a certain diffusion process with ``multiplication''
jumps. The paper also includes three ``elementary'' proofs for the
simple prototype equation $y'(x)+y(x)=\frac{1}{2}\mypp y(q\myp
x)+\frac{1}{2}\mypp y(x/q)$, based on perturbation, analytical, and
probabilistic techniques, respectively, which may appear useful in
other situations as efficient exploratory tools.


\end{abstract}

% 34K05 General theory
% 34K06 Linear functional-differential equations
% 34K12 Growth, boundedness, comparison of solutions
% 45J05 Integro-ordinary differential equations
% 60Jxx Markov processes
% 60J25 Markov processes with continuous parameter
% 60J35 Transition functions, generators and resolvents [See also 47D03, 47D07]
% 60G44 Martingales with continuous parameter
% 60G46 Martingales and classical analysis
% 60G40 Stopping times; optimal stopping problems; gambling theory [See also 62L15, 91A60]
% 39A13 Difference equations, scaling ($q$-differences)
% 39-xx:Difference and functional equations
% 34E20 Singular perturbations, turning point theory, WKB methods

\section{Introduction}\label{sec0}
The classical \emph{pantograph equation} is the linear first-order
functional-differential equation (with rescaled argument) of the
form
\begin{equation}\label{eq0:pantograph}
y'(x)=a\myp y(q\myp x)+b\myp y(x),
\end{equation}
where $a,b$ are constant coefficients (real or complex) and $q>0$ is
a rescaling parameter. Historically,\footnote{The name ``pantograph
equation'' was not in wide use until it was coined by Iserles
\cite{I} for a more general class of functional-differential
equations.} the term ``pantograph'' dates back to the seminal paper
of 1971 by Ockendon and Tayler \cite{OT}, where such
equations\footnote{To be more precise, a certain vector analog of
Eq.~(\ref{eq0:pantograph}).} emerged in a mathematical model for the
dynamics of an overhead current collection system on an electric
locomotive (with the physically relevant value $q <1$). At about the
same time, a systematic analysis of solutions to the pantograph
equation was started by Fox \textit{et al}.\ \cite{FMOT}, where
various analytical, perturbation, and numerical techniques were
discussed at length (for both $q <1$ and $q >1$).

Interestingly, an equation of the form (\ref{eq0:pantograph}) (with
$q >1$) was derived more than 25 years earlier by Ambartsumian
\cite{Amb} to describe the absorption of light by the interstellar
matter. Some particular cases of Eq.~(\ref{eq0:pantograph}) are also
found in early work by Mahler \cite{Mah} on a certain partition
problem in number theory (where Eq.~(\ref{eq0:pantograph}), with
$a=1$, $b=0$, $q <1$, appears as a limit of a similar
functional-difference equation) and by Gaver \cite{Gav} on a special
ruin problem (with $a=1$, $b=-1$, $q >1$). Subsequently, the
pantograph equation has appeared in numerous applications ranging
from the problem of coherent states in quantum theory \cite{Sp} to
cell-growth modeling in biology \cite{Wake} (see further references
in Refs.\ \cite{Der1,DM,I,Marsh2}). These and other examples suggest
that, typically, the pantograph equation and similar
functional-differential equations with rescaling are relevant as
long as the systems in question possess some kind of
self-similarity.

Since its introduction into the mathematical literature in the early
1970s, the theory of the pantograph equation (and some of its
natural generalizations) has been the subject of persistent
attention and research effort, yielding over years a number of
significant developments. In particular, the classification of
Eq.~(\ref{eq0:pantograph}) with regard to various domains of the
parameters,\footnote{Depending on whether $q <1$ or $q >1$ and also
on the cases $\Re b <0$, $\Re b >0$, and $\Re b =0$.} including
existence and uniqueness theorems, and an extensive asymptotic
analysis of the corresponding solutions have been given by Kato and
McLeod \cite{KM} and Kato \cite{K}. The investigation of such
equations in the complex domain was initiated by Morris \emph{et
al.}\ \cite{MFB} and Oberg \cite{Oberg} and continued by Derfel and
Iserles \cite{DI} and Marshall \emph{et al.}\ \cite{Marsh2}. A
systematic treatment of the generalized first-order pantograph
equation (with matrix coefficients and also allowing for a term with
rescaled derivative) is contained in the influential paper by
Iserles \cite{I}, where in particular a fine geometric structure of
almost-periodic solutions has been described. Asymptotics for
equations with variable coefficients have been studied by Derfel and
Vogl \cite{DV}.

Higher-order generalizations of the pantograph equation
(\ref{eq0:pantograph}) lead to the class of linear
functional-differential equations with rescaling,
\begin{equation}\label{eq0:fde}
\sum_{j=1}^\ell \sum_{k=0}^m a_{jk}\mypp y^{(k)}(\alpha_j\myp
x+\beta_j)=0
\end{equation}
(see Ref.\ \cite{Der2} and further references
therein).\footnote{Note that theory of such equations is closely
related to the theory of $q$-difference equations developed by
Birkhoff \cite{B} and Adams \cite{A} (see also Section \ref{sec3}
below).} Kato \cite{K} posed a problem of asymptotic analysis of
Eq.~(\ref{eq0:fde}), including the question of existence and
characterization of \emph{bounded} solutions. Some partial answers
to the latter question have been given by Derfel \cite{Der1,Der2}
and Derfel and Molchanov \cite{DM}.

In particular, Derfel \cite{Der1} considered the ``balanced''
generalized first-order pantograph equation of the form
\begin{equation}\label{eq0:alpha}
y'(x)+y(x)=\sum_{j=1}^\ell p_j\mypp y(\alpha_j\myp x),
\end{equation}
subject to the condition
\begin{equation}\label{eq:balance}
\sum_{j=1}^\ell p_j=1,\ \ \quad p_j> 0 \quad(j=1,\dots,\ell),
\end{equation}
so that the weights $p_j$ of the rescaled $y$-terms on the
right-hand side of Eq.~(\ref{eq0:alpha}) match the unit coefficient
of the $y(x)$ on the left. Note that, owing to the balance condition
(\ref{eq:balance}), Eq.~(\ref{eq0:alpha}) always has a trivial
solution $y=\const$. The question of existence of \emph{nontrivial}
(i.e., nonconstant) \emph{bounded solutions} is most interesting
(and most difficult) in the case where the right-hand side of
Eq.~(\ref{eq0:alpha}) involves both ``advanced'' ($\alpha_j>1$) and
``delayed'' ($0<\alpha_j<1$) arguments. It turns out that the answer
depends crucially on the quantity
\begin{equation}\label{eq:K}
K:=\sum_{j=1}^\ell p_j\ln \alpha_j\,.
\end{equation}
Namely, Derfel \cite{Der1} has proved that if $K<0$ then
Eq.~(\ref{eq0:alpha}) has no nontrivial bounded solutions, whereas
if $K>0$ then such a solution always exists. In the ``critical''
case $K=0$, this question has remained open as yet.

In the present paper, we consider a more general
integro-differential equation\footnote{In fact, the results of the
paper \cite{Der1} mentioned above include first-order equations of
the form (\ref{eq:int-dif}), i.e., with $a_2=0$. Let us also remark
that more general first-order integro-differential equations (but
with delayed arguments only, i.e., $\alpha\in(0,1)$) were considered
by Iserles and Liu \cite{IL}.} of the pantograph type, namely,
\begin{equation}\label{eq:int-dif}
-a_2\myp y''(x)+a_1 y'(x)+y(x)=\int_0^\infty y(\alpha x)\,\mu(\rd
\alpha),
\end{equation}
where $a_1\ge0$, \,$a_2\ge0$, \,$a_1^2+a_2^2>0$ (so that $a_1,a_2$
do not vanish simultaneously), and $\mu$ is a probability measure on
$(0,\infty)$,
\begin{equation}\label{eq:mu-bal}
\mu(0,\infty)=\int_0^\infty \mu(\rd\alpha)=1.
\end{equation}
The parameter $\myp\alpha\myp$ in Eq.~(\ref{eq:int-dif}) can be
viewed as a random variable, with values in $(0,\infty)$ and the
probability distribution given by the measure $\mu$, i.e.,
$\PP\{\alpha\in A\}=\mu(A)$, $A\subset(0,\infty)$. Note that
Eq.~(\ref{eq:int-dif}) is balanced in the same sense as
Eq.~(\ref{eq0:alpha}), since the mean contribution of the
distributed rescaled term $y(\alpha x)$ is matched by that of
$y(x)$. Moreover, Eq.~(\ref{eq:int-dif}) reduces to
Eq.~(\ref{eq0:alpha}) when $a_1=1$, $a_2=0$, and the measure $\mu$
is discrete, with atoms $p_j=\mu(\alpha_j)$, $j=1,\dots,\ell$
\,(i.e., $\alpha$ is a discrete random variable, with the
distribution $\PP\{\alpha=\alpha_j\}=p_j$, $j=1,\dots,\ell$).

As already mentioned, due to the balance condition (\ref{eq:mu-bal})
any constant satisfies Eq.~(\ref{eq:int-dif}), and by linearity of
the equation one can assume, without loss of generality, that
$y(0)=0$. Moreover, if $x>0$ ($x<0$) then the right-hand side of
Eq.~(\ref{eq:int-dif}) is determined solely by the values of the
function $y(u)$ with $u>0$ (respectively, $u<0$). Therefore, the
two-sided equation (\ref{eq:int-dif}) is decoupled at $x=0$ into two
one-sided boundary value problems,
\begin{equation}\label{eq:bv}
\begin{aligned} -&a_2\myp y''(x)+a_1 y'(x)+y(x)=\int_0^\infty y(\alpha
x)\,\mu(\rd \alpha),\qquad
x\gtrless 0,\\
&y(0)=0.
\end{aligned}
\end{equation}

For Eq.~(\ref{eq:int-dif}), the analog of Eq.~(\ref{eq:K}) is given
by
\begin{equation}\label{eq:K1}
K:=\int_0^\infty \ln
\alpha\:\mu(\rd\alpha)=\EE\myp[\myn\ln\alpha\myp].
\end{equation}
Our main result is the following theorem, which resolves the problem
of nontrivial bounded solutions in the critical case, $K=0$ (and
also recovers and extends the result of Ref.\ \cite{Der1} for the
case $K<0$).
\begin{theorem}\label{th1.1}
Assume that $0\ne \EE\myp|\ln\alpha|<\infty$, so that $K$ in
Eq.~(\ref{eq:K1}) is well defined and the measure $\mu$ is not
concentrated at the point\/ $\alpha=1$, i.e., the random variable
$\alpha$ does not degenerate to the constant~$1$. Under these
hypotheses, the condition $K\le 0$ implies that any bounded solution
of equation (\ref{eq:int-dif}) is trivial, i.e., $y(x)\equiv\const$,
$x\in\RR$.
\end{theorem}

The apparent probabilistic structure of Eq.~(\ref{eq:int-dif}) is
crucial for our proof of this result. The main idea is to construct
a certain diffusion process $X_t$, with negative drift and
``multiplication'' jumps (i.e., of the form $x\mapsto\alpha x$),
such that Eq.~(\ref{eq:int-dif}) can be rewritten as $\calL y=0$,
where $\calL$ is the infinitesimal generator of the Markov process
$X_t$. That is to say, the class of bounded solutions of
Eq.~(\ref{eq:int-dif}) coincides with the set of bounded
$\calL$-harmonic functions. This link brings in the powerful tool
kit of Markov processes; particularly instrumental is the well-known
fact (see, e.g., Ref.\ \cite{EK}) that for any $\calL$-harmonic
function $f(x)$, the random process $f(X_t)$ is a martingale, and
hence, for any $t\ge0$,
\begin{equation}\label{eq:mart00}
f(x)=\EE\myp[f(X_t)\myp|X_0=x],\qquad x\in\RR.
\end{equation}
On the other hand, due to the multiplication structure of
independent consecutive jumps of the process $X_t$, its position
after $n$ jumps is expressed in terms of a background random walk
$S_k=\xi_1+\cdots+\xi_k$ ($0\le k\le n$), where $\xi_i$'s are
independent random variables with the same distribution as
$\ln\alpha$. The hypothesis $K\le 0$ of Theorem \ref{th1.1} implies
that, almost surely (a.s.), the random walk $S_n$ travels
arbitrarily far to the left. Using an optional stopping theorem
(whereby the boundedness of $f(x)$ is important), we can apply the
martingale identity (\ref{eq:mart00}) at the suitably chosen
stopping (first-passage) times, which eventually leads to the
conclusion that $f'(x)\equiv0$ and hence $f(x)\equiv\const$.

This approach also allows us to give an example of a nontrivial
bounded solution to equation (\ref{eq:int-dif}) in the case $K>0$
(thus extending the result by Derfel \cite{Der1} to the second-order
pantograph equation).
\begin{theorem}\label{th1.2}
Suppose that $K>0$, and set
\begin{equation}\label{f-infty}
f_\infty(x):=\PP\Bigl\{\liminf_{t\to\infty} X_t=+\infty\myp|\myp
X_0=x\Bigr\},\qquad x\in\RR,
\end{equation}
where
$X_t$ is the random process constructed in the proof of
Theorem~\ref{th1.1}. Then the function $f_\infty(x)$ is
$\calL$-harmonic and such that $f_\infty(x)\to0$ as $x\to-\infty$
and $f_\infty(x)\to1$ as $x\to+\infty$.
\end{theorem}

In the particular case where $a_2=0$, Eq.~(\ref{eq:int-dif}) becomes
\begin{equation}\label{eq:int-dif1}
a_1 y'(x)+y(x)=\int_0^\infty y(\alpha x)\,\mu(\rd \alpha),\qquad
x\in\RR.
\end{equation}
In this case, the diffusion component of the random process $X_t$ is
switched off, and it follows, due to the negative drift and
multiplication jumps (see details in Section \ref{sec4}), that if
$X_0=x\le0$, then $X_t\le 0$ for all $t\ge0$. That is to say, the
negative semi-axis $(-\infty,0]$ is an absorbing set for the process
$X_t$, and hence the function $f_\infty(x)$, defined by
Eq.~(\ref{f-infty}) as the probability to escape to $+\infty$
starting from $x$, vanishes for all $x\le 0$. This leads to the
following interesting specification of the example in Theorem
\ref{th1.2}.
\begin{corollary}
If\/ $a_2=0$ in Eq.~(\ref{eq:int-dif}) then $f_\infty(x)\equiv 0$
for all $x\le 0$. Moreover, Eq.\ (\ref{eq:int-dif}) implies that all
derivatives  of the function $f_{\infty}(x)$ vanish at zero,
$f^{(k)}_{\infty}(0)=0$ $(k=1,2,\dots)$.
\end{corollary}

Before elaborating the ideas outlined above, we would like to make a
short digression in order to consider the simple prototype example
of Eq.~(\ref{eq0:alpha}), namely,
\begin{equation}\label{eq0:1/2}
y'(x)+y(x)=\frac{1}{2}\,y(q\myp x)
+\frac{1}{2}\,y(q^{-1}x)\qquad(q\ne1),
\end{equation}
and to give several different ``sketch'' proofs of Theorem
\ref{th1.1} in this case. Note that, according to Eq.~(\ref{eq:K}),
we have
$$
K=\frac{1}{2}\ln q+\frac{1}{2}\ln\frac{1}{q}=0,
$$
so Eq.~(\ref{eq0:1/2}) falls in the (most interesting) critical
case. In fact, this example was the starting point of our work and a
kind of mathematical test-tube to try various approaches and ideas.
Although not strictly necessary for the exposition, after some
deliberation we have cautiously decided to include our early proofs
(based on perturbation, analytical, and probabilistic
arguments,\footnote{It is amusing that these three methods represent
nicely the traditional organization of British mathematics into
Applied Mathematics, Pure Mathematics, and Statistics, which is
reflected in the names of mathematical departments in most
universities in the UK.} respectively), partly because this will
hopefully equip the reader with some insight into validity of the
result, and also because these methods may appear useful as
exploratory tools in other situations.

The rest of the paper is laid out as follows. In Sections
\ref{sec2}, \ref{sec3}, and \ref{sec4}, we discuss the three
approaches to equation (\ref{eq0:1/2}) as just mentioned. In Section
\ref{sec5} we start a more systematic treatment by describing the
construction of a suitable diffusion process with multiplication
jumps. In Section \ref{sec6}, we discuss the corresponding
$\calL$-harmonic functions and obtain an \emph{a priori} bound for
the derivative of a solution. Finally, in Section \ref{sec7} we
prove our main Theorems \ref{th1.1} and~\ref{th1.2}.

\section{Perturbative proof}\label{sec2}

Following the ideas used by Ockendon and Tayler \cite{OT} and Fox
\emph{et al.}\ \cite{FMOT} in the case of the original pantograph
equation (\ref{eq0:pantograph}), we start by observing that if $q=1$
then Eq.~(\ref{eq0:1/2}) is reduced to the equation $y'=0$, which
has constant solutions only. Therefore, when the parameter $q$ is
close to $1$, it is reasonable to seek solutions of
Eq.~(\ref{eq0:1/2}) in a form that involves ``superposition'' of
(exponentially) small oscillations (\emph{fast variation}) on top of
an almost constant (polynomial) function (\emph{slow variation}).
This leads to a WKB\myp-type asymptotic expansion of the solution in
terms of perturbation parameter $\varepsilon\approx0$ (see Ref.\
\cite{Shiva}), which in the first-order approximation yields two
first-order differential equations: a nonlinear equation (called the
\emph{eikonal equation}) for the fast variation  and a linear
equation (called the \emph{transport equation}) for the slow
variation. For simplicity of presentation, we will restrict
ourselves to the first-order approximation, but in principle one can
go on to the analysis of higher-order terms, which are described by
linear equations and therefore can be determined without much
trouble.

To implement this approach, set $q=1\pm\varepsilon$
\,($\varepsilon>0$), \,$x=\varepsilon^{-1}u$, and
$y(x)=y(\varepsilon^{-1}u)=:f(u)$. Then
$$
y'(x)=\frac{\rd f}{\rd u}\cdot \frac{\rd u}{\rd x}=\varepsilon
f'(u),
$$
and Eq.~(\ref{eq0:1/2}) takes the form
\begin{equation}\label{eq:q1}
\varepsilon f'(u)+f(u)=\frac{1}{2}\,f\bigl((1\pm\varepsilon)\myp u
\bigr)+\frac{1}{2}\,f\bigl((1\pm\varepsilon)^{-1}u\bigr).
\end{equation}
As explained in the Introduction, without loss of generality we may
assume that $f(0)=0$.

Now, suppose that, for small $\varepsilon$, the function $f(u)$
admits a WKB\myp-type expansion,
\begin{equation}\label{eq:q0}
f(u)\sim \bigl(A_0(u)+\varepsilon A_1(u)+\cdots\bigr)
\myp\exp\bigl(\varepsilon^{-1}V(u)\bigr),
\end{equation}
which in principle should be valid uniformly for all $u$, including
the limiting values $u\to0$ and $u\to\infty$. Differentiation of
Eq.~(\ref{eq:q0}) yields
\begin{equation}\label{eq:eps}
\varepsilon f'(u)\sim \bigl(A_0(u) V'(u)+\varepsilon
\bigl(A'_0(u)+A_1(u)V'(u)\bigr)+\cdots\bigr)\myp\exp\bigl(\varepsilon^{-1}V(u)\bigr).
\end{equation}
From Eq.~(\ref{eq:q0}) we also obtain
\begin{equation}\label{eq:e+}
\begin{aligned}
f\bigl((1\pm\varepsilon)\myp u\bigr) &\sim \bigl(A_0(u)+ \varepsilon
\bigl(A_1(u)\pm uA'_0(u)\bigr)+\cdots\bigr)\\
&\quad\times \exp\left(\frac{V(u)}{\varepsilon}\pm
uV'(u)+\frac{\varepsilon u^2 V''(u)}{2}+\cdots\right),
\end{aligned}
\end{equation}
and
\begin{equation} \label{eq:e-}
\begin{aligned}
&f\bigl((1\pm\varepsilon)^{-1} u\bigr)\sim \bigl(A_0(u)+
\varepsilon \bigl(A_1(u)\mp uA'_0(u)\bigr)+\cdots\bigr)\\
&\qquad\quad\times \exp\left(\frac{V(u)}{\varepsilon}\mp
uV'(u)+\varepsilon
\left(uV'(u)+\frac{u^2\,V''(u)}{2}\right)+\cdots\right).
\end{aligned}
\end{equation}
Substituting the expansions (\ref{eq:eps}), (\ref{eq:e+}) and
(\ref{eq:e-}) into Eq.~(\ref{eq:q1}), canceling out the common
factor $\exp\bigl(V(u)/\varepsilon\bigr)$, and collecting the terms
that remain after setting $\varepsilon=0$, we get
$$
A_0(u) V'(u) + A_0(u)=\frac{1}{2}\,A_0(u)\exp\bigl(\pm uV'(u)\bigr)
+\frac{1}{2}\,A_0(u)\exp\bigl(\mp uV'(u)\bigr).
$$
Assuming that $A_0(u)\ne 0$, this gives the equation
\begin{equation}\label{eq:q2}
1+V'(u)=\cosh\bigl(uV'(u)\bigr),
\end{equation}
or equivalently
\begin{equation}\label{eq:q3}
u=\frac{uV'(u)}{\cosh\bigl(uV'(u)\bigr)-1\mathstrut}.
\end{equation}
Similarly, equating the terms of order of $\varepsilon$ and noting
that $A_1(u)$ cancels out owing to Eq.~(\ref{eq:q2}), we obtain
\begin{align}\label{eq:q4}
A'_0(u)=A'_0(u)\,u\sinh\bigl(uV'(u)\bigr)&+\frac{1}{2}\,A_0(u)\,
u^2\,V''(u)\cosh \bigl(uV'(u)\bigr)\\
\notag
&+\frac{1}{2}\,A_0(u)\, u V'(u)\exp\bigl(\mp uV'(u)\bigr)
\end{align}

We can now check that the formal expansion (\ref{eq:q0}) is
compatible with the zero initial condition, $f(0)=0$. Equation
(\ref{eq:q3}) implies that if $u\to0$ then $uV'(u)\to\infty$, and
moreover
\begin{equation}\label{eq:qq1}
uV'(u)\sim -\ln u+\ln\ln\frac{1}{u}+\cdots,
\end{equation}
whence
\begin{equation}\label{eq:V}
V(u)\sim -\frac{\ln^2\! u}{2}+\ln u \cdot\ln \ln\frac{1}{u}+\cdots.
\end{equation}
Furthermore, differentiation of Eq.~(\ref{eq:qq1}) gives
\begin{equation}\label{eq:qq3}
\begin{aligned}
u^2V''(u) &\sim \ln u-\ln\ln\frac{1}{u} +\cdots.
\end{aligned}
\end{equation}
Inserting formulas (\ref{eq:qq1}) and (\ref{eq:qq3}) into
Eq.~(\ref{eq:q4}), we obtain for $A_0(u)$ the asymptotic
differential equation
$$
\frac{A_0'(u)}{A_0(u)}\sim\mp\frac{1}{2u}\left(\ln u
-\ln\ln\frac{1}{u}+\cdots\right),
$$
which solves to
\begin{equation}\label{eq:A0}
\ln A_0(u)\sim \mp\frac{\ln u}{4}\left(\ln u- 2
\ln\ln\frac{1}{u}+\cdots\right).
\end{equation}
Finally, substituting the expansions (\ref{eq:V}) and (\ref{eq:A0})
into Eq.~(\ref{eq:q0}) we obtain that $f(u)\to 0$ as $u\to0$, as
required.

Let us now explore the behavior of the solution as $u\to\infty$. In
this limit, Eq.~(\ref{eq:q3}) gives $uV'(u)\to0$, and moreover
$$
u\sim \frac{2}{uV'(u)}-\frac{uV'(u)}{6}+\cdots,
$$
whence
\begin{align}
\label{eq:V'}
u V'(u)&\sim \frac{2}{u}-\frac{2}{3u^3}+\cdots,\\
\label{eq:V''}
u^2V''(u)&\sim -\frac{4}{u}+\frac{8}{3u^3}+\cdots.
\end{align}
From (\ref{eq:V'}), we also find
\begin{equation}\label{eq:V1}
V(u)\sim C-\frac{2}{u}+\frac{2}{9u^3}+\cdots,
\end{equation}
where $C=\const$. Inserting the expansions (\ref{eq:V'}) and
(\ref{eq:V''}) into Eq.~(\ref{eq:q4}), we obtain
\begin{align*}
\frac{A'_0(u)}{A_0(u)}\sim\frac{1}{u}\pm \frac{2}{u^2}+\cdots,
\end{align*}
and hence
\begin{equation}\label{eq:Cu}
A_0(u)\sim C_0\myp u\exp\left(\mp\frac{2}{u}+ \cdots\right).
\end{equation}
The case $C_0\ne0$ is unsuitable, since Eq.~(\ref{eq:Cu}) would
imply that $A_0(u)\to\infty$ as $u\to\infty$ and, in view of
formulas (\ref{eq:V1}) and (\ref{eq:q0}), the solution $f(u)$
appears to be unbounded, which contradicts our assumption.
Therefore, $C_0=0$ and hence $A_0(u)=0$, thus reducing  the
expansion (\ref{eq:q0}) to
\begin{equation*}
f(u)\sim \varepsilon \bigl(A_1(u)+\varepsilon A_2(u)+\cdots\bigr)
\myp\exp\bigl(\varepsilon^{-1}V(u)\bigr).
\end{equation*}
Arguing as above, we successively obtain $A_1(u)=0$, $A_2(u)=0$,
etc. This indicates that $f(u)=0$, which was our aim.

\section{Analytical proof}\label{sec3}

In this section, we demonstrate how the theory of $q$-difference
equations (see Refs.\ \cite{A,B}) can be used to show that
Eq.~(\ref{eq0:1/2}) has no nontrivial bounded solutions. In what
follows, we assume that $q\ne1$.

As explained in the Introduction (see Eq.~(\ref{eq:bv})),
Eq.~(\ref{eq0:1/2}) splits into two (similar) one-sided equations,
so it suffices to consider the boundary-value problem
\begin{equation}\label{eq0:00}
\begin{aligned}
&y'(x)+y(x)=\frac{1}{2}\,y(q\myp x) +\frac{1}{2}\,y(q^{-1}x),\qquad x\ge0,\\
&y(0)=0.
\end{aligned}
\end{equation} Assume that $y(x)$ is a
bounded solution of Eq.~(\ref{eq0:00}), $|y(x)|\le B$ $(x\ge 0)$,
and let $\hat y(s)$ be the Laplace transform of $y(x)$,
$$
\hat y(s):=\int_{0}^\infty \re^{-s x}\myp y(x)\,\rd x,
$$
then $\hat y(s)$ is analytic in the right half-plane, $\Re s>0$, and
\begin{equation}\label{eq:Re}
|\hat y(s)|\le \frac{B}{\Re s},\qquad \Re s>0.
\end{equation}
On account of the boundary condition $y(0)=0$, Eq.~(\ref{eq0:00})
transforms into
\begin{equation}\label{eq:fq3'}
(1+s)\myp \hat y(s)=\frac{1}{2\myp q}\,\hat
y(q^{-1}s)+\frac{q}{2}\,\hat y(qs),
\end{equation}
or, after the substitution $\varphi(s):=s\mypp \hat y(s)$,
\begin{equation}\label{eq:fq3}
(1+s)\myp
\varphi(s)=\frac{1}{2}\mypp\varphi(q^{-1}s)+\frac{1}{2}\mypp\varphi(qs).
\end{equation}
Note that the estimate (\ref{eq:Re}) implies
\begin{equation}\label{eq:Re1}
|\varphi(s)|\le \frac{B\myp |s|}{\Re s},\qquad \Re s>0,
\end{equation}
and in particular $\varphi(s)$ is bounded in the vicinity of the
origin.

Let us rewrite Eq.~(\ref{eq:fq3}) in the form
\begin{equation}\label{eq:fq4}
\varphi(q^2s)-2(1+qs)\myp \varphi(qs)+\varphi(s)=0.
\end{equation}
Equation (\ref{eq:fq4}) is a linear $q$-difference equation of
order~$2$. According to the general theory of such equations (see
Ref.\ \cite{A}), the characteristic equation for Eq.~(\ref{eq:fq4})
(in the vicinity of $s=0$) reads
\begin{equation*}
q^{2\rho}-2q^{\myp\rho}+1=0,
\end{equation*}
and $\rho_{1,2}=0$ is its multiple root. The corresponding
fundamental set of solutions to Eq.~(\ref{eq:fq4}) is given by
\begin{align*}
\varphi_{1}(s)&=P_1(s),\\
\varphi_{2}(s)&=P_2(s)+\frac{\ln s}{\ln q}\,P_3(s),
\end{align*}
where $P_1(s)$, $P_2(s)$, and $P_3(s)$ are generic power series
convergent in some neighborhood of zero, and $\ln s$ denotes the
principal branch of the logarithm. Note that the solution
$\varphi_{2}(s)$ is unsuitable because it is unbounded near $s=0$
(see Eq.~(\ref{eq:Re1})). On the other hand, the function
$\varphi_{1}(s)$ is analytic in the vicinity of zero and, moreover,
it can be analytically continued, step by step, into the whole
complex plane $\CC$ by means of Eq.~(\ref{eq:fq3}). For instance, if
$q>1$ then the analytic continuation from a disk $|s|\le a$ to the
bigger disk $|s|\le q\myp a$ is furnished by the formula
$\varphi(qs)=2\myp(1+s)\mypp \varphi(s)-\varphi(q^{-1}s)$ (see
Eq.~(\ref{eq:fq3})), and so on.

That is to say, $\varphi_{1}(s)$ can be extended to an entire
function $\varphi(s)$, which by construction satisfies
Eq.~(\ref{eq:fq4}) for all $s\in\CC$. But then, according to one
result by Mason \cite{Mason}, the entire function $\varphi(s)$ must
be of \emph{zero order}, and consequently (see Ref.\
\cite[\S\myp8.7.3]{Titch}) it is \emph{unbounded} on any ray (in
particular, for $s\in(0,\infty)$), unless it is a constant. However,
the unboundedness for real $s>0$ contradicts the estimate
(\ref{eq:Re1}). Hence, $\varphi(s)=\const$, so that $\hat
y(s)=\const\cdot\myp s^{-1}$, and by the uniqueness theorem for the
Laplace transform this implies that $y(x)\equiv\const$, i.e,
$y(x)\equiv y(0)=0$, as claimed.

\section{Probabilistic proof}\label{sec4}

In this section, we give a probabilistic interpretation of
Eq.~(\ref{eq0:1/2}) via a certain ruin problem, and prove that the
corresponding solution is constant using elementary probabilistic
considerations. Although our argument does not cover the whole class
of bounded solutions, it contains some ideas that we will use in the
second half of the paper to give a complete proof of our general
result.

Let us consider the following ``double-or-half'' gambling model (in
continuous time). Suppose that a player spends his initial capital,
$x$, at rate $v$ per unit time, so that after time $t$ he is left
with capital $x-vt$. However, at a random time $\tau$ (with
exponential distribution), he gambles by putting the remaining
capital at stake, whereby he can either \emph{double} his money or
lose \emph{half} of it, both with probability $1/2$. After that, the
process continues in a similar fashion, independently of the past
history. If the capital reaches zero and then moves down to become
negative, this is interpreted as borrowing, so the process proceeds
in the same way without termination. In that case, gambling will
either double or halve the debt, and in particular the capital will
remain negative forever.

More generally, if $X_t$ denotes the player's capital at time
$t\ge0$, starting with the initial amount $X_0=x$, then the random
process $X_t$ moves with constant negative drift $(-v)$, interrupted
at random time instants $\sigma_i$ by random multiplication jumps
from its location $x_i=X_{\sigma_i-0}$ (i.e., immediately before the
jump) to either $q\myp x_i$ or $q^{-1}x_i$ ($q\ne1$), both with
probability $1/2$. We assume that the jumps occur at the arrival
times $\sigma_1, \sigma_2,\dots$ of an auxiliary Poisson process
with parameter $\lambda>0$, so that the waiting times until the next
jump, $\tau_i=\sigma_i-\sigma_{i-1}$ ($\sigma_0:=0$), are
independent identically distributed (i.i.d.) random variables, each
with the exponential distribution
$$
\PP\{\tau>t\}=\re^{-\lambda t},\qquad t>0.
$$
According to this description, $(X_t)$ is a Markov process, in that
the probability law of its future development is completely
determined by its current state, but not by the past history (``lack
of memory'') (see, e.g., Refs.\ \cite{EK,Fe}).

We are concerned with the \emph{ruin problem} for this
model.\footnote{A similar ruin problem for the process with
deterministic multiplication jumps of the form $x_i\mapsto q\myp
x_i$ ($q>1$), was first considered by Gaver \cite{Gav}, leading to
the equation $y'(x)+y(x)=y(qx)$ (i.e., with advanced argument, cf.\
Eq.~(\ref{eq0:alpha})). The systematic theory of general processes
with multiplication jumps was developed by Lev \cite{Lev}.} Namely,
consider the probability $f_0(x)$ of becomong bankrupt starting with
the initial capital $x$,
$$
f_{0}(x):=\PP\Bigl\{\liminf_{t\to\infty} X_t\le 0\,|\mypp
X_0=x\Bigr\}\equiv\PP_x\{T_0<\infty\},\qquad x\in\RR,
$$
where $T_0:=\min\{t\ge 0: X_t\le 0\}$ is the random time to
bankruptcy and $\PP_x$ denotes the probability measure conditioned
on the initial state $X_0=x$.

From the definition of the process $X_t$, it is clear that if $x\le
0$ then $T_0=0$ and so $f_0(x)=1$. For $x>0$, we note that if the
first jump does not occur prior to time $x/v$, then the process will
simply drift down to $0$, in which case $T_0=x/v<\infty$. Otherwise
(i.e., if a jump does happen before time $x/v$), the ruin problem
may be reformulated by treating the landing point after the jump as
a new starting point (thanks to the Markov property). More
precisely, by conditioning on the first jump instant $\sigma_1$
($=\tau_1$) and using the (strong) Markov property, we obtain
\begin{align}
\notag
f_{0}(x)&=\PP_x\{\sigma_1>x/v\}+\int_0^{x/v}
\lambda\mypp\re^{-\lambda s}\,\PP_x(T_0<\infty\,|\,\sigma_1=s)\,\rd s\\
\notag &= \re^{-\lambda x/v}+\lambda\int_0^{x/v} \re^{-\lambda
s}\,\frac{1}{2}\mypp\biggl(f_0\bigl(q(x-vs)\bigr)+
f_0\bigl(q^{-1}(x-vs)\bigr)\biggr)\,\rd s\\
\label{eq:repr}  &=\re^{-\lambda
x/v}\left(1+\frac{\lambda}{2v}\int_0^x \re^{\lambda
u/v}\bigl(f_0(qu)+f_0(q^{-1}u) \bigr)\,\rd u \right),
\end{align}
where in the last line we have made the substitution $u=x-vs$. The
representation (\ref{eq:repr}) implies that the function $f_{0}(x)$
is continuous and, moreover, (infinitely) differentiable, and by
differentiation of Eq.~(\ref{eq:repr}) with respect to $x$, it
follows that the function $y=f_{0}(x)$ satisfies the generalized
pantograph equation (cf.\ Eq.~(\ref{eq0:1/2}))
\begin{equation}\label{eq:fdeq}
\frac{v}{\lambda}\mypp y'(x)+y(x)=\frac{1}{2}\,y(q\myp
x)+\frac{1}{2}\,y(q^{-1}x).
\end{equation}
It is easy to see that, in fact, this equation is satisfied on the
whole axis, $x\in\RR$. As we have mentioned, $f_0(x)=1$ for all
$x\le 0$, and it is now our aim to show that the same is true for
all $x>0$, which would mean that the solution $y=f_0(x)$ to equation
(\ref{eq0:1/2}) is a constant, $f_0(x)\equiv 1$, $x\in\RR$.

To this end, note that the position of the process after $n$ jumps
is given by
$$
X_{\sigma_n}=(X_{\sigma_{n-1}}-v\myp\tau_n)\,q^{\myp\xi_n},\qquad
n=1,2,\dots,
$$
where $\xi_n$'s are i.i.d.\ random variables taking the values
$\pm1$ with probabilities $1/2$. By iterations (using that $X_0=x$),
we obtain\footnote{Similar random sums as in Eq.~(\ref{eq:Z}) arise
in products of certain random matrices in relation to random walks
on the group of affine transformations of the line (see Ref.\
\cite{Gr}).}
\begin{align}
\notag
X_{\sigma_n}&=\Bigl(\bigl((x-v\myp\tau_1)\,q^{\myp\xi_1}-v\myp\tau_2\bigr)\mypp
q^{\myp\xi_2}-\dots-
v\myp\tau_n\Bigr)\mypp q^{\myp\xi_n}\\
\label{eq:Z} &=(x-v\myp\tau_1)\,q^{\myp\xi_1+\xi_2+\dots+\xi_n}-
v\myp\tau_2\mypp q^{\myp\xi_2+\dots+\xi_n}-\dots-
v\myp\tau_n\myp q^{\myp\xi_n}\\
\notag &=q^{\myp S_n} \biggl(x-v\sum_{i=1}^n\tau_i\mypp
q^{-S_{i-1}}\biggr),
\end{align}
where $S_n:=\xi_1+\xi_2+\dots+\xi_n$, \,$S_0:=0$. Note that $S_n$
can be interpreted as a (simple) random walk, which in our case is
symmetric (i.e., $\PP\{\xi_i=1\}=\PP\{\xi_i=-1\}=1/2$) and therefore
\emph{recurrent} (see, e.g., Ref.\ \cite{Fe}). In particular, the
events $A_n:=\{S_{n-1}=0\}$ ($n=1,2,\dots$) occur infinitely often,
with probability~$1$. Furthermore, setting $B_n:=\{\tau_n>1\}$, we
note that the events $A_n\cap B_n$ ($n=1,2,\dots$) are conditionally
independent, given the realization of the random walk $\{S_k,\, k\ge
1\}$. Since the random variables $\tau_n$ (and therefore the events
$B_n$) are independent of $\{S_k\}$, we have, with probability $1$,
\begin{align*}
\sum_{n=1}^\infty \PP(A_n\cap B_n\myp|\myp \{S_k\})&=
\sum_{n=1}^\infty \PP(B_n)\,{\bf
1}_{A_n}\\&=\re^{-\lambda}\myp\#\{n:A_n\ \text{occurs}\}=\infty,
\end{align*}
where ${\bf 1}_{\{\cdots\}}$ denotes the indicator of an event.
Hence, Borel-Cantelli's lemma (see, e.g., Ref.\ \cite{Fe}) implies
\begin{equation}\label{eq:a.s.}
\PP(A_n\cap B_n\ \text{occur infinitely often}\mypp|\myp
\{S_k\})=1\qquad\text{(a.s.)},
\end{equation}
and by taking the expectation in Eq.~(\ref{eq:a.s.}) (with respect
to the distribution of the sequence $\{S_k\}$), the same is true in
the unconditional form,
$$
\PP(A_n\cap B_n\ \text{occur infinitely often})=1.
$$
As a consequence, the terms in the random series
$$
\sum_{n=1}^\infty \tau_n\myp q^{-S_{n-1}}
$$
will infinitely often exceed the value $1$, all other terms being
nonnegative. Therefore, the series diverges to $+\infty$ (a.s.), and
from (\ref{eq:Z}) it follows that
$$
\liminf_{n\to\infty} X_{\sigma_n}\le 0\qquad (\text{a.s.}).
$$
In turn, this implies that $T_0<\infty$ (a.s.), and so $f_0(x)=1$
for all $x>0$, as claimed.

\section{Jump diffusions}\label{sec5}

We now pursue a more general (and more systematic) approach.
Equations of the form (\ref{eq0:fde}) are linked in a natural way
with certain continuous-time Markov processes (more specifically,
diffusions with multiplication jumps). To describe this class of
processes, let us consider a Brownian motion $B^{\kappa,v}_t$,
starting at the origin, with diffusion coefficient $\kappa\ge 0$ and
nonpositive (constant) drift $-v\le0$,
$$
\rd B^{\kappa,v}_t=\kappa\myp {\rd B}_t-v\myp \rd t,\qquad
B^{\kappa,v}_0=0,
$$
or equivalently
$$
B^{\kappa,v}_t=\kappa B_t-v\myp t, \qquad t\ge0,
$$
where $B_t=B_t^{1,\myp0}$ is a standard Brownian motion (with
continuous sample paths). We assume that $\kappa^2+v^2>0$, so that
$B^{\kappa,v}_t$ does not degenerate to a (zero) constant.

The random process $B_t^{\kappa,v}$ determines the underlying
diffusion dynamics for a process with jumps, $(X_t,\,t\ge0)$, which
is defined as follows. Suppose that the jump instants are given by
the arrival times $\sigma_1$, $\sigma_2,\dots$ of an auxiliary
Poisson process with parameter $\lambda>0$, so that
$\tau_i=\sigma_{i}-\sigma_{i-1}$ \,($i=1,2,\dots$) are i.i.d.\
random variables with exponential distribution,
$$
\PP\{\tau>t\}=\re^{-\lambda t},\qquad t>0
$$
(we set formally $\sigma_0:=0$). Furthermore, suppose that the
successive jumps are determined by the rescaling coefficients
$\alpha_i$ of the form $\alpha_i=\re^{\myp\xi_i}$, where $\xi_i$'s
are i.i.d.\ random variables. Then, the (right-continuous) sample
paths of the process $X_t$ are defined inductively by
\begin{equation}\label{eq:D}
X_t=\left\{ \begin{array}{ll} x+B^{\kappa,v}_t,&\ \ 0=\sigma_0\le t<\sigma_{1},\\[.3pc]
\re^{\myp\xi_i}
X_{\sigma_i-0}+B^{\kappa,v}_{t}-B^{\kappa,v}_{\sigma_i},&\ \
\sigma_{i}\le t<\sigma_{i+1},\ \ \quad i=1,2,\dots
\end{array}
\right.
\end{equation}
That is to say, the process $(X_t,\,t\ge0)$ starts at point $x$ and
moves as $X_t=x+B_t^{\kappa,v}$ until a random time $\sigma_1$, when
it jumps to a random point
$$
X_{\sigma_1}= \re^{\myp\xi_1} X_{\sigma_1-0}= \re^{\myp\xi_1}
(x+B_{\sigma_1}^{\kappa,v})=\re^{\myp\xi_1}
x+\re^{\myp\xi_1}\zeta_1,
$$
where
$\zeta_1:=B_{\sigma_1}^{\kappa,v}=B_{\sigma_1}^{\kappa,v}-B_{\sigma_0}^{\kappa,v}$.
Thereafter, the process proceeds in a diffusive way as
$X_t=X_{\sigma_1}+B^{\kappa,v}_{t}-B^{\kappa,v}_{\sigma_1}$ until a
random time $\sigma_2$, when it makes the next jump to
\begin{align*}
X_{\sigma_2}&=\re^{\myp\xi_2}X_{\sigma_2-0}=\re^{\myp\xi_2}
(X_{\sigma_1}+B^{\kappa,v}_{\sigma_2}-B^{\kappa,v}_{\sigma_1})\\
&=\re^{\myp\xi_2+\xi_1}
x+\re^{\myp\xi_2+\xi_1}\zeta_1+\re^{\myp\xi_2}\zeta_2,
\end{align*}
where $\zeta_2:=B_{\sigma_2}^{\kappa,v}-B_{\sigma_{1}}^{\kappa,v}$,
and so on. Iterating, we obtain that the $n$-th jump, occurring at
time $\sigma_n$, lands at the point
\begin{align}
\label{eq:X_tau} X_{\sigma_n}&
=\re^{\myp\xi_n}X_{\sigma_n-0}=\re^{\myp\xi_n}
(X_{\sigma_{n-1}}+B^{\kappa,v}_{\sigma_n}-B^{\kappa,v}_{\sigma_{n-1}})\\
\notag & =\re^{\myp S_n}(x+\zeta_1+\zeta_2\mypp\re^{-S_1}+\cdots
+\zeta_n\mypp\re^{-S_{n-1}}),
\end{align}
where
\begin{equation}\label{eq:zeta}
\begin{aligned}
\zeta_n:={}&B_{\sigma_n}^{\kappa,v}-B_{\sigma_{n-1}}^{\kappa,v}\\
{}={}&\kappa\myp(B_{\sigma_n}-B_{\sigma_{n-1}})-v\myp(\sigma_n-\sigma_{n-1})
\qquad (n=1,2,\dots)
\end{aligned}
\end{equation}
is a sequence of i.i.d.\ random variables.

\section{$\calL$-harmonic functions}\label{sec6}
Let us study the properties of the process $X_t$ in greater detail.
\begin{definition}
\normalfont The \emph{infinitesimal operator} (\emph{generator})
$\calL$ of a Mar\-kov random process $(X_t,\,t\ge 0)$ is defined by
\begin{equation}\label{eq:L0}
(\calL f)(x):=\lim_{h\to0+}\frac{\EE_x[f(X_h)]-f(x)}{h},
\end{equation}
with the domain $\calD(\calL)$ consisting of functions $f$ for which
the limit in Eq.~(\ref{eq:L0}) exists.
\end{definition}

In a standard way, by considering possible scenarios for the process
$(X_t)$ up to an infinitesimal time $h$ (see Ref.\ \cite{EK}), one
obtains the following.

\begin{proposition}
For the random diffusion with jumps $(X_t)$ defined above, its
generator $\calL$ acts on bounded\/ $C^2$-smooth functions as
\begin{equation}\label{eq:L}
(\calL f)(x)=\frac{\kappa^2}{2}\myp f''(x)-vf'(x)+\lambda\bigl(\EE
[f(\re^{\myp\xi} x)]- f(x)\bigr),
\end{equation}
where $\xi$ is a random variable with the same distribution as any
one of the i.i.d.\ random variables $\xi_1,\xi_2,\dots$.
\end{proposition}

\begin{definition}
\normalfont A function $f$ is called \emph{$\calL$-harmonic} if\/
$\calL f=0$.
\end{definition}

Note that the expectation in Eq.~(\ref{eq:L}) can be written as a
Stieltjes integral,
$$
\EE [f(\re^{\myp\xi} x)]=\int_0^\infty \!f(z x)\,\rd F(z),
$$
where $F(z)$ is the cumulative distribution function of the random
variable $Z=\re^{\myp\xi}$, i.e., $F(z):=\PP\{\re^{\myp\xi}\le
z\}=\PP\{\xi\le\ln z\}$ \,($0<z<\infty$). Hence, the  equation
$\calL f=0$, with $\calL$ given by Eq.~(\ref{eq:L}), is equivalent
to
$$
-\frac{\kappa^2}{2}\myp f''(x)+vf'(x)+\lambda
f(x)=\lambda\int_0^\infty \!f(z x)\,\rd F(z),
$$
which is a balanced generalized pantograph equation of the form
(\ref{eq:int-dif}).

Our aim is to study the class of bounded $\calL$-harmonic functions.
Let us denote by $\|\cdot\|$ the standard sup-norm on $\RR$:
$$
\|f\|:=\sup_{x\in\RR} |f(x)|.
$$
As a first step, we estimate the derivative of an $\calL$-harmonic
function.

\begin{proposition}\label{pr:f'} Suppose that $v>0$.
If\/ $\|f\|<\infty$ and $\calL f=0$ then $\|f'\|<\infty$.
\end{proposition}

\proof If $\kappa=0$ then, according to Eq.~(\ref{eq:L}), the
equation $\calL f=0$ takes the form
$$
vf'(x)=\lambda\myp\bigl(\EE\myp[f(\re^{\myp\xi} x)]-f(x)\bigr),
$$
which gives
$$
\|f'\|\le \frac{2\lambda}{|v|}\,\|f\|<\infty.
$$

For $\kappa>0$, the condition $\calL f=0$ is equivalent to
\begin{equation}\label{eq:f''}
f''(x)-\gamma f'(x)=-g(x),
\end{equation}
where
$$
\gamma:=\frac{2v}{\kappa^2}>0,\qquad
g(x):=\frac{2\lambda}{\kappa^2}\,\bigl(\EE\myp[f(\re^{\myp\xi}x)]
-f(x)\bigr).
$$
Solving equation (\ref{eq:f''}), we obtain
\begin{equation}\label{eq:f'}
\begin{aligned}
f'(x)&=\int_x^\infty g(u)\, \re^{-\gamma (u-x)}\,\rd u
+C\myp\re^{\gamma
x}\\
&=\int_0^\infty\! g(u+x)\, \re^{-\gamma u}\,\rd u+C\myp\re^{\gamma
x},
\end{aligned}
\end{equation}
where $C=\const$. Since $\|g\|\le 4\lambda \kappa^{-2}\|f\|<\infty$,
Eq.~(\ref{eq:f'}) implies that if $C\ne0$ then
$$
f'(x)=O(1)+C\myp\re^{\gamma x}\to\infty \qquad(x\to+\infty),
$$
so that $\lim_{x\to+\infty}f(x)=\infty$, which contradicts the
assumption $\|f\|<\infty$. Therefore, $C=0$ and
$$
\|f'\|\le \|g\|\int_0^\infty\! \re^{-\gamma u}\,\rd
u=\frac{\|g\|}{\gamma}<\infty.
$$
The proof is complete. \endproof

\section{Proof of the main results}\label{sec7}
Let $(\calF_{t},\,t\ge0)$ be the natural filtration generated by the
process $(X_t)$, i.e., $\calF_t=\sigma\{X_s,\,s\le t\}$ is the
minimal $\sigma$-algebra containing all ``level'' events $\{X_s\le
c\}$ ($c\in\RR$, $s\le t$). Intuitively, $\calF_t$ is interpreted as
the collection of all the information that can be obtained by
observation of the random process $(X_s)$ up to time $t$. As is well
known (see, e.g., Ref.\ \cite[Ch.~4]{EK}), if a function $f$ is
$\calL$-harmonic then the random process $f(X_t)$ is a
\emph{martingale} relative to $(\calF_t)$, i.e., for any $0\le s\le
t$, with probability~$1$,
\begin{equation}\label{eq:mart}
\EE_{x}[f(X_t)\myp|\mypp\calF_s]=f(X_s).
\end{equation}
In words, this means that if we are trying to predict the mean value
of the martingale $f(X_t)$ at some future time $t$ using its past
history up to time $s\le t$, then the best estimate is given by
$f(X_s)$, i.e., the value of the process at the latest available
time instant $s$.

In particular, taking expectation of both sides of
Eq.~(\ref{eq:mart}) at $s=0$ gives
\begin{equation}\label{eq:mart0}
\EE_{x}[f(X_t)]=f(x),\qquad t\ge0.
\end{equation}
In addition, if $f$ is bounded then, by Doob's optional stopping
theorem (see, e.g., Ref.\ \cite[\S\myp8]{Yeh}), Eq.~(\ref{eq:mart0})
extends to
\begin{equation}\label{eq:marT}
\EE_{x}[f(X_T)]=f(x),
\end{equation}
where $T$ is a random \emph{stopping time}, i.e., such that $\{T\le
t\}\in\calF_t$ for each $t\ge0$. In words, one should be able to
decide whether the random time $T$ has occurred by observing the
process up to a given time. The martingale property (\ref{eq:marT})
is crucial in the proof of our main theorem below, where we will
apply it to the special sequence of stopping times.

\begin{theorem}[cf.\ Theorem \ref{th1.1}]\label{th:const}
Assume that\/ $\EE\myp [\myp\xi\myp]\le0$ and $\PP\{\xi\ne 0\}>0$.
If\/ $\calL f=0$ and\/ $\|f\|<\infty$, then $f(x)\equiv \const$,
$x\in\RR$.
\end{theorem}
\proof
 For \,$r>0$, set
$N_r:=\min\{n:S_n\le -r\}$, where $S_n=\xi_1+\dots +\xi_n$, and
define
\begin{equation}\label{eq:sigma'}
T_r:=\sigma_{N_r}=\sum_{n=1}^\infty \sigma_n\mypp{\bf
1}_{\{N_r=n\}},
\end{equation}
where $(\sigma_i)$ are the time instants of successive jumps (see
Section~\ref{sec3}). The assumption $\EE\myp [\myp\xi\myp]\le 0$
implies (see, e.g., Ref.\ \cite{Fe}) that $\liminf_{n\to\infty}
S_n=-\infty$ (a.s.), so that $N_r<\infty$ (a.s.) and therefore the
random variable $T_r$ is well defined.

Note that $T_r$ is a stopping time for the random process $(X_t)$,
i.e., $\{{T_r\le t}\}\allowbreak\in\calF_t$ for each $t\ge0$.
Indeed, suppose that we are given a sample path of the process $X_s$
up to time $t$, and in particular we know the time instants
$\sigma_i$ and the magnitudes of all the jumps prior to $t$. Then,
using Eq.~(\ref{eq:X_tau}), we can reconstruct the corresponding
values\footnote{There is a slight problem if
$X_{\sigma_i-0}=X_{\sigma_i}=0$, but this only happens with zero
probability, so may be ignored.} $\xi_i=\ln
(X_{\sigma_i}/X_{\sigma_i-0})$. In turn, this allows us to determine
if the threshold $(-r)$ has been reached by the associated random
walk $S_n$ and, therefore, whether or not the condition $T_r\le t$
holds, as required.

Now, applying the optional stopping theorem in the form
(\ref{eq:marT}), we obtain
\begin{equation}\label{eq:f}
f(x)=\EE_x [f(X_{T_r})]=\EE\bigl[f(\re^{S_{N_r}} x+\text{terms
independent of $x$})\bigr].
\end{equation}
First, suppose that $v>0$. Differentiation of Eq.~(\ref{eq:f}) with
respect to $x$ gives
\begin{equation}\label{eq:f'1}
f'(x)=\EE\bigl[\myp\re^{S_{N_r}} f'(\re^{S_{N_r}} x+\cdots)\bigr],
\end{equation}
whence, using Proposition \ref{pr:f'}, we get
$$
\|f'\|\le \re^{-r}\myp\|f'\|<\infty.
$$
Letting here $r\to+\infty$, we conclude that $\|f'\|=0$ and so
$f=\const$.

If $v=0$ then, differentiating Eq.~(\ref{eq:f'1}) we get
$$
f''(x)=\EE\bigl[\myp\re^{2S_{N_r}} f''(\re^{S_{N_r}}
x+\cdots)\bigr],
$$
so that
\begin{equation}\label{eq:f''1}
\|f''\|\le \re^{-2r}\|f''\|<\infty.
\end{equation}
Equation (\ref{eq:f''}) (with $\gamma=0$) implies
$\|f''\|\le\|g\|<\infty$, so taking $r\to+\infty$ in
Eq.~(\ref{eq:f''1}) yields $f''(x)\equiv 0$. Therefore,
$f(x)=c_1x+c_0$, but since $\|f\|<\infty$, we must have $c_1=0$, so
that $f(x)\equiv c_0=\const$.
\endproof

\begin{proposition}\label{pr:f}
Denote $A_\infty:=\bigl\{\liminf_{t\to\infty} X_t=+\infty\bigr\}$,
and consider the probability of the event $A_\infty$ as a function
of the initial point of the process $X_t$,
$$
f_{\infty}(x):=\PP_{x}(A_\infty)=\PP(A_\infty\myp|\myp X_0=x),\ \
\quad x\in \RR.
$$
Then the function $f_{\infty}(x)$ is $\calL$-harmonic, i.e.,
$$
(\calL f_{\infty})(x)=0,\ \ \quad x\in \RR,
$$
where $\calL$ is the generator of the random process $(X_t)$ given
by Eq.~(\ref{eq:L}).
\end{proposition}
\proof Conditioning on $X_h$ ($<+\infty$), by the Markov property we
obtain
\begin{align*}
f_{\infty}(x)
&=\EE_x[\PP_x(A_\infty\myp|\myp X_h)]\\
&=\EE_x [\PP_{X_h}(A_\infty)]\\
&=\EE_x [f_\infty(X_h)],
\end{align*}
whence, by the definition (\ref{eq:L0}), it readily follows that
$\calL f_\infty=0$.
\endproof

\begin{theorem}[cf.\ Theorem \ref{th1.2}]
Suppose that $\EE\myp[\myp\xi\myp]>0$. Then the function
$f_\infty(x)$ is a nontrivial bounded $\calL$-harmonic function; in
particular, $f_\infty(x)\to 0$ as $x\to-\infty$ and $f_\infty(x)\to
1$ as $x\to+\infty$.
\end{theorem}

\proof  The function $f_\infty(x)$ is $\calL$-harmonic by
Proposition \ref{pr:f}. In order to obtain the limits of
$f_\infty(x)$ as $x\to\pm\infty$, note that
\begin{align}
\label{eq:i.o.} 1-f_\infty(x)&=\PP_{x}\Bigl\{\liminf_{t\to\infty}
X_t<+\infty\Bigr\}\\
\notag
&=\PP_{x}\Bigl\{\liminf_{n\to\infty} X_{\sigma_n}<+\infty\Bigr\}\\
\notag&= \lim_{M\to\infty}\PP_{x} \{X_{\sigma_n}\le M\
\text{infinitely often}\}.
\end{align}
According to Eq.~(\ref{eq:X_tau}), the condition $X_{\sigma_n}\le M$
can be rewritten as
\begin{equation}\label{eq:M}
x+\zeta_1+\zeta_2\mypp\re^{-S_1}+\cdots
+\zeta_n\mypp\re^{-S_{n-1}}\le M\mypp\re^{-S_n}.
\end{equation}
Since $\EE\myp[\myp\xi\myp]>0$, the strong Law of Large Numbers
implies that, with probability~$1$,
\begin{equation}\label{eq:LLN}
S_n\sim n\myp\EE\myp[\myp\xi\myp]\to+\infty\qquad (n\to\infty).
\end{equation}
It follows that if the inequality (\ref{eq:M}) holds for infinitely
many $n$, then
\begin{equation}\label{eq:liminf}
x+\liminf_{n\to\infty}\left(\zeta_1+\zeta_2\mypp\re^{-S_1}+\cdots
+\zeta_n\mypp\re^{-S_{n-1}}\right)\le \lim_{n\to\infty}
M\mypp\re^{-S_n}=0.
\end{equation}
Moreover, using Eq.~(\ref{eq:LLN}) and recalling that $\zeta_i$ are
i.i.d.\ random variables (see Eq.~(\ref{eq:zeta})), it is easy to
show (e.g., using Kolmogorov's ``three series'' theorem, see Ref.\
\cite{Fe}) that the random series
\begin{equation}\label{eq:Upsilon}
\eta:=\sum_{n=1}^\infty \zeta_n\mypp\re^{-S_{n-1}}
\end{equation}
converges with probability~$1$. Therefore, from Eqs.\ (\ref{eq:M}),
(\ref{eq:liminf}) and (\ref{eq:Upsilon}) it follows that for any
$M>0$,
$$
\PP_{x}\{X_{\sigma_n}\le M\ \text{infinitely often}\}\le
\PP\{\eta\le -x\}.
$$
Returning to Eq.~(\ref{eq:i.o.}), we deduce that
\begin{align*}
1\ge f_\infty(x)&\ge 1-\PP\{\eta\le -x\}\\
&=\PP\{\eta>-x\}\to 1\qquad (x\to+\infty).
\end{align*}

On the other hand, writing the left-hand side of Eq.~(\ref{eq:M}) as
$x+\eta+\delta_n$, where $\delta_n\to0$ (a.s.), we have, for any
$\varepsilon>0$, $M>0$,
\begin{align*}
\PP\{x+\eta\le -\varepsilon\}&\le \PP\{x+\eta+\delta_n\le
0\ \,\text{for all $n$ large enough}\}\\
&\le \PP\{x+\eta+\delta_n< M\mypp\re^{-S_n}\ \text{infinitely
often}\},
\end{align*}
which in view of Eq.~(\ref{eq:i.o.}) implies
\begin{align*}
0\le f_\infty(x)&\le 1-\PP\{\eta+x\le
-\varepsilon\}\\
&=\PP\{\eta>-x-\varepsilon\}\to 0\qquad (x\to-\infty).
\end{align*}
Thus, the proof is complete.
\endproof

\section*{Acknowledgments}

%The second author (G.D.) wishes to thank Arieh Iserles for
%stimulating discussions during and useful remarks. The third author
%(S.M.) gratefully acknowledges the support from the Center of
%Advanced Studies in Mathematics of the Ben Gurion University during
%his visit in May--July 2006, and he is thankful to Michael Lin for
%the kind hospitality.

Part of this research was done when the second author (G.D.) was
visiting the University of Cambridge in May--June 2005, and his
thanks are due to Arieh Iserles for stimulating discussions and
useful remarks. The third author (S.M.) gratefully acknowledges the
support from the Center of Advanced Studies in Mathematics of the
Ben Gurion University during his visit in May--June 2006, and he
would like to thank Michael Lin for kind hospitality.

% Editor: please remove the next line when not needed
\pagebreak


\begin{thebibliography}{99}
\bibitem {A}
{\sc C.R.~Adams}, {\sl Linear $q$-difference equations}, Bull.\
Amer.\ Math.\ Soc.\ {\bf 37} (1931), 361--400.

\bibitem{Amb}
{\sc V.A.~Ambartsumian}, {\sl On the theory of brightness
fluctuations in the Milky Way}, (Russian) Doklady Akad.\ Nauk SSSR
{\bf 44} (1944), 244--247; (English translation) Compt.\ Rend.\
(Doklady) Acad.\ Sci.\ URSS {\bf 44} (1944), 223--226.

\bibitem {B}
{\sc G.D.~Birkhoff},  {\sl The generalized Riemann problem for
linear differential equations and the allied problems for linear
difference and $q$-difference equations}, Proc.\ Amer.\ Acad. Arts
Sci.\ {\bf 49} (1913), 521--568.

\bibitem {Der1}
{\sc G.A.~Derfel}, {\sl Probabilistic method for a class of
functional-differential equations\-}, (Russian) Ukrain.\ Mat.\ Zh.\
{\bf 41} (1989),
%No. 10,
1322--1327; (English translation) Ukrainian\ Math.\ J.\ {\bf 41}
(1990), 1137--1141.

\bibitem {Der2}
%{\sc G.~Derfel},
{\sc G.~Derfel}, {\sl Functional-differential and functional
equations with rescaling}, in Operator Theory and Boundary
Eigenvalue Problems (International Workshop, Vienna, July 27-30,
1993), Operator Theory: Advances and Applications, Vol.~80,
\,I.~Gohberg and H.~Langer (eds.), Birkh\"auser, Basel, 1995, pp.\
100--111.

\bibitem {DI}
{\sc
%\bysame\
G.A.~Derfel and A.~Iserles}, {\sl The pantograph equation in the
complex plane}, J.~Math.\ Anal.\ Appl.\ {\bf 213} (1997), 117--132.

\bibitem {DM}
{\sc
%\bysame\
G.~A.~Derfel and S.A.~Molchanov}, {\sl Spectral methods in the
theory of func\-tional-differential equations}, (Russian) Mat.\
Zametki {\bf 47}\myp(3) (1990),
%no. 3,
42--51; (English translation) Math.\ Notes {\bf 47} (1990),
254--260.
%no. 3-4,

\bibitem {DV}
{\sc
%\bysame\
G.~Derfel and F.~Vogl}, {\sl On the asymptotics of solutions of a
class of linear  functional-differen\-tial equations}, European J.\
Appl.\ Math.\ {\bf 7} (1996), 511--518.
%no. 5

\bibitem{EK}
{\sc S.N.~Ethier and T.G.~Kurtz}, {\sl Markov Processes:
Characterization and Convergence}, Wiley Series in Probability and
Mathematical Statistics, John Wiley \& Sons, New York, 1986.

\bibitem{Fe}
{\sc W.~Feller}, {\sl An Introduction to Probability Theory and Its
Applications}, Vol.~II, Wiley Series in Probability and Mathematical
Statistics, 2nd ed., John Wiley \& Sons, New York, 1971.

\bibitem{FMOT}
{\sc L.~Fox, D.F.~Mayers, J.R.~Ockendon, and A.B.~Tayler}, {\sl On a
functional dif\-ferential equation}, J.~Inst.\ Math.\ Appl.\ {\bf 8}
(1971), 271--307.

\bibitem{Gav}
{\sc D.P.\ Gaver, Jr.},  {\sl An absorption probability problem},
J.~Math.\ Anal.\ Appl.\ {\bf 9} (1964), 384--393.
% no. 3

\bibitem{Gr}
{\sc A.K.~Grintsevichyus}, {\sl  On the continuity of the
distribution of a sum of dependent variables connected with
independent walks on lines}. (Russian) Teor.\ Veroyatn.\ i
Primenen.\ {\bf 19} (1974), 163--168; (English translation\ Teor.\
Probab.\ Appl.\ {\bf 19} (1974), 163--168.

\bibitem{I} {\sc A.~Iserles}, {\sl On the generalized pantograph
functional-differential equation}, European J.\ Appl.\ Math.\ {\bf
4} (1993), 1--38.
%no. 1

\bibitem{IL}
{\sc
%\bysame\
A.~Iserles and Y.K.~Liu}, {\sl On pantograph integro-differential
equations}, J.\ Integral Equations Appl.\ {\bf 6} (1994),
%no. 2,
213--237.

\bibitem{K}
{\sc T.~Kato}, {\sl Asymptotic behavior of solutions of the
functional differential equation $y'(x)=a\myp y(\lambda x)+b\myp
y(x)$}, in Delay and Functional Differential Equations and Their
Applications (Proc.\ Conf., Park City, Utah, March 6--11, 1972),
K.~Schmitt (ed.), Academic Press, New York, 1972, pp.\ 197--217.

\bibitem{KM}
{\sc T.~Kato and J.B.~McLeod}, {\sl The functional-differential
equation $y'(x)=a\myp y(\lambda x)+b\myp y(x)$}, Bull.\ Amer.\
Math.\ Soc.\ {\bf 77} (1971), 891--937.

\bibitem{Lev}
{\sc G.Sh.~Lev}, {\sl Semi-Markov processes of multiplication with
drift},  (Russian) Teor.\ Veroyatn.\ i Primenen.\ {\bf 17} (1972),
160--166; (English translation) Theory Probab.\ Appl.\ {\bf 17}
(1972), 159--164.
%no. 3-4,

\bibitem{Mah}
{\sc K.~Mahler}, {\sl On a special functional equation}, J.~London
Math.\ Soc.\ {\bf 15} (1940), 115--123.

\bibitem{Marsh2}
{\sc J.C.~Marshall, B.~van-Brunt, and G.C.~Wake}, {\sl  A natural
boundary for solutions to the second order pantograph equation},
J.~Math.\ Anal.\ Appl.\ {\bf 299} (2004),
%no. 2,
314--321.

\bibitem{Mason}
{\sc T.E.~Mason}, {\sl On properties of the solutions of linear
$q$-difference equations with entire function coefficients}, Amer.\
J.\ Math.\ {\bf 37} (1915), 439--444.

\bibitem{MFB}
{\sc G.R.~Morris, A.~Feldstein, and E.W.~Bowen}, {\sl The
Phragm\'en-Lindel\"of principle and a class of functional
differential equations}, in Ordinary Differential Equations (Proc.\
Conf., Math.\ Res.\ Center, Naval Res.\ Lab., Washington, D.C.,
1971), L.\,Weiss (ed.), Academic Press, New York, 1972, pp.\
513--540.

\bibitem{Oberg}
{\sc R.J.~Oberg}, {\sl Local theory of complex functional
differential equations}, Trans.\ Amer.\ Math.\ Soc.\ {\bf 161}
(1971), 302--327.

\bibitem{OT} {\sc J.R.~Ockendon and A.B.~Tayler}, {\sl The
dynamics of a current collection system for an electric locomotive},
Proc.\ Royal Soc.\ London A {\bf 322} (1971), 447--468.

\bibitem{Shiva}
{\sc B.K.~Shivamoggi}, {\sl Perturbation Methods for Differential
Equations}, Birkh\"auser, Boston, 2003.

\bibitem{Sp}
{\sc V.~Spiridonov}, {\sl Universal superpositions of coherent
states and self-similar potentials}, Phys.\ Rev.\ A {\bf 52} (1995),
1909--1935.

\bibitem{Titch}
{\sc E.C.~Titchmarsh}, {\sl The Theory of Functions}, 2nd ed.,
Oxford University Press, Oxford, 1939.

\bibitem{Wake}
{\sc G.C.~Wake, S.~Cooper, H.K.~Kim, and B.~van-Brunt}, {\sl
Functional differential equations for cell-growth models with
dispersion}, Commun.\ Appl.\ Anal.\ {\bf 4} (2000), 561--573.

\bibitem{Yeh}
{\sc J.~Yeh}, {\sl Martingales and Stochastic Analysis}, Series on
Multivariate Analysis, Vol.\,1, World Scientific, Singapore, 1995.

\end{thebibliography}
\end{document}